\def\yes{\if00}

\def\iftwelvept{\yes}
\def\ifusepdf{\yes}
\def\ifusepsfont{\yes}

\iftwelvept
\documentclass[leqno,12pt]{amsart}
\else
\documentclass[leqno]{amsart}
\fi
\usepackage{amssymb}
\usepackage{amscd}
\usepackage{latexsym}
\usepackage{verbatim}
\ifusepdf
\usepackage{hyperref}
\else\fi
\ifusepsfont
\usepackage[T1]{fontenc}
\usepackage{times}
\else\fi

\iftwelvept
\else\fi

\theoremstyle{plain}
\newtheorem{Theorem}{Theorem}[section]

\newtheorem{Lemma}[Theorem]{Lemma}

\newtheorem{Claim}{Claim}[Theorem]
\newtheorem{Formula}[Claim]{Formula}

\theoremstyle{definition}

\renewcommand{\theTheorem}{\arabic{section}.\arabic{Theorem}}

\renewcommand{\theequation}{\arabic{section}.\arabic{Theorem}.\arabic{Claim}}

\def\rom{\textup}

\newcommand{\QQ}{{\mathbb{Q}}}

\newcommand{\CC}{{\mathbb{C}}}
\newcommand{\PP}{{\mathbb{P}}}
\newcommand{\OO}{{\mathcal{O}}}
\newcommand{\XX}{{\mathcal{X}}}
\newcommand{\LL}{{\mathcal{L}}}
\newcommand{\Pic}{\operatorname{Pic}}
\newcommand{\aPic}{\widehat{\operatorname{Pic}}}

\newcommand{\zero}{\operatorname{div}}
\newcommand{\Proof}{{\sl Proof.}\quad}
\newcommand{\adeg}{\widehat{\operatorname{deg}}}
\newcommand{\trdeg}{\operatorname{tr.deg}}
\newcommand{\acherncl}{\widehat{{c}}}
\newcommand{\normabb}{\Vert\!\cdot\!\Vert}
\newcommand{\QED}{{\unskip\nobreak\hfil\penalty50\quad\null\nobreak\hfil
{$\Box$}\parfillskip0pt\finalhyphendemerits0\par\medskip}}
\newcommand{\rest}[2]{\left.{#1}\right\vert_{{#2}}}
\begin{document}

%%%%%%%%%%%
%% Title %%
%%%%%%%%%%%
\title[The continuity of Deligne's pairing]%
{The continuity of Deligne's pairing}
\author{Atsushi Moriwaki}
\address{Department of Mathematics, Faculty of Science,
Kyoto University, Kyoto, 606-8502, Japan}
\email{moriwaki@kusm.kyoto-u.ac.jp}
\date{29/April/1999, 10:10PM (JP), (Version 2.01)}
%\keywords{}
%\subjclass{}

\maketitle

\section*{Introduction}
\renewcommand{\theTheorem}{\Alph{Theorem}}

In the paper \cite[\S1.2]{ZhHRSV},
S. Zhang expected the metric of Deligne's pairing to be continuous.
In this note, we will give an affirmative answer for his question.
Namely, we will show the following theorem.

\begin{Theorem}
\label{thm:cont:metric:Deligne:pairing}
Let $f : X \to S$ be a flat and projective morphism of
algebraic varieties over $\CC$ with $n = \dim f$.
Let $\overline{L}_0, \ldots, \overline{L}_n$ be $C^{\infty}$-hermitian
line bundles on $X$. Then, the metric of
Deligne's pairing $\langle \overline{L}_0, \ldots, \overline{L}_n
\rangle(X/S)$ is continuous.
\end{Theorem}

As an application of the above theorem,
we have the following result, which is, in some sense,
a generalization of Kawaguchi's Hodge index theorem \cite{KaHI}.

\begin{Theorem}
\label{thm:hodge:index}
Let $K$ be a finitely generated field over $\QQ$ with
$d = \trdeg_{\QQ}(K)$, and let
$\overline{B} = (B, \overline{H})$ be a polarization of $K$
\rom{(}for details, see \cite{MoArht}\rom{)}.
Let $X$ be a smooth projective
curve of genus $g$ over $K$,
$J_X$ the Jacobian of $X$, and let $\Theta_X$ be a symmetric theta divisor on $J_X$.
For a line bundle $L$ on $X$ with $\deg(L) = 0$,
there is a sequence of $C^{\infty}$-models $(\XX_m, \overline{\LL}_m)$
of $(X, L)$ over $B$ with 
\[
\lim_{n\to\infty} \adeg \left(
\acherncl_1 (\overline{\mathcal{L}}_m) \cdot
\acherncl_1 (\overline{\mathcal{L}}_m) \cdot
\acherncl_1(f_m^*(\overline{H}))^d \right) = 
-2 \hat{h}^{\overline{B}}_{\Theta_X}([L]),
\]
where $f_m$ is the canonical morphism $\XX_m \to B$.
\end{Theorem}

Finally, we would like to express hearty thanks to Prof. Kawaguchi for telling us
the exact formula \eqref{eqn:step2:proof:main:theorem}.
The author also thanks Prof. Zhang for his nice comments.

\section{Deligne's pairing with metric}
\renewcommand{\theTheorem}{\arabic{section}.\arabic{Theorem}}

Here we give a quick review of Deligne's pairing with metric.
For details, see \cite[\S8.1]{D}, \cite{E1}, \cite{E2}, \cite{F1}, \cite{F2},
\cite[\S1.1-\S1.2]{ZhHRSV} and etc.

Let $f : X \to S$ be a flat and projective morphism of
integral schemes of relative dimension $n$.
Then, we have Deligne's pairing
\[
\langle \ , \ldots, \ \rangle(X/S) :
 \overbrace{\Pic(X) \times \cdots \times \Pic(X)}^{\text{$(n+1)$-times}}
 \to \Pic(S).
\]
Roughly speaking, Deligne's pairing is given by
\[
 c_1\left(\langle L_0, \ldots, L_n \rangle(X/S) \right) = 
f_*(c_1(L_0) \cdots c_1(L_n))
\]
for $L_0, \ldots, L_n \in \Pic(X)$.
Note that
if $L_0, \ldots, L_n$ have local (with respect to $S$)
sections $l_0, \ldots, l_n$ without intersections,
then the symbol $\langle l_0, \ldots, l_n \rangle$ 
gives rise to a local base of
$\langle L_0, \ldots, L_n \rangle(X/S)$.

Next we assume that $X$ and $S$ are defined over $\CC$.
Let $\overline{L}_0, \ldots, \overline{L}_n$ be 
$C^{\infty}$-hermitian line bundles on $X$.
The metric of
$\langle L_0, \ldots, L_n \rangle(X/S)$ induced
by metrics of $\overline{L}_0, \ldots, \overline{L}_n$
is defined inductively in the following way.
Let $l_0, \ldots, l_n$ be local 
sections of $L_0, \ldots, L_n$ such that
$l_0, \ldots, l_n$ have no intersections and
$Y = \zero(l_n)$ is integral and flat over $S$.
Then, 
the length of $\langle l_0, \ldots, l_n \rangle$
is given by
\[
 \Vert \langle \rest{l_0}{Y}, \ldots, \rest{l_{n-1}}{Y} \rangle \Vert
 \exp\left( \int_{X/S} \log \Vert l_n \Vert
 \bigwedge_{i=0}^{n-1} c_1(\overline{L}_i )\right),
\]
where
\[
 \left( \int_{X/S} \log \Vert l_n \Vert
 \bigwedge_{i=0}^{n-1} c_1(\overline{L}_i )\right)(s)
 =
 \int_{X_s} \log \Vert l_n \Vert
 \bigwedge_{i=0}^{n-1} c_1(\overline{L}_i )
\]
for each $s \in S(\CC)$.
We denote the line bundle
$\langle L_0, \ldots, L_n \rangle(X/S)$ with the above metric by
$\langle \overline{L}_0, \ldots, \overline{L}_n \rangle(X/S)$.
For simplicity,
\[
 \langle 
 \overbrace{\overline{L}_0, \ldots, \overline{L}_0}^{\text{$a_0$-times}},
 \overbrace{\overline{L}_1, \ldots, \overline{L}_1}^{\text{$a_1$-times}},
 \ldots, 
 \overbrace{\overline{L}_t, \ldots, \overline{L}_t}^{\text{$a_t$-times}} 
 \rangle(X/S)
\]
is denoted by $\langle \overline{L}_0^{\cdot a_0}, \overline{L}_1^{\cdot a_1}, \ldots,
\overline{L}_t^{\cdot a_t} \rangle(X/S)$.

\section{The continuity of the fiber integral}
\label{sec:cont:fiber:integral}
In this section, we will consider the continuity of
the fiber integral of $C^{\infty}$-forms.
First of all, let us fix the $C^{\infty}$ of hermitian line
bundles in this note, which coincides with
the sense of \cite{ZhHRSV} and \cite{MoArht}.

Let $X$ be an algebraic variety over $\CC$, and $\overline{L}$ a
continuous hermitian line bundle on $X$.
We say $\overline{L}$ is a {\em $C^{\infty}$-hermitian line bundle} if,
for any complex manifolds $M$ and any analytic maps $f : M \to X$,
$f^*(\overline{L})$ is a $C^{\infty}$-hermitian line bundle on $M$.
In the same way, a continuous function $\phi$ on $X$ is said to be
{\em $C^{\infty}$} if,
for any complex manifolds $M$ and any analytic maps $f : M \to X$,
$f^*(\phi)$ is a $C^{\infty}$-function on $M$.

The proof of the following theorem is the main purpose of this section.

\begin{Theorem}
\label{thm:cont:fiber:integral}
Let $f : X \to S$ be a flat and projective morphism of
algebraic varieties over $\CC$ with $n = \dim f$.
Let $\overline{L}_1, \ldots, \overline{L}_n$ be $C^{\infty}$-hermitian
line bundles on $X$, and let $\phi$ be a $C^{\infty}$-function on $X$. 
Then, the fiber integral
${\displaystyle
 \int_{X/S} \phi \bigwedge_{i=1}^n c_1(\overline{L}_i)
}$
is continuous, namely, the function given by
\[
S(\CC) \ni s \mapsto  \int_{X_s} \phi \bigwedge_{i=1}^n c_1(\overline{L}_i)
\]
is continuous.
\end{Theorem}

\Proof
Before starting the proof of our theorem,
we would like to introduce the $C^{\infty}$ of
hermitian line bundles in strong sense.
A continuous function $\phi$ on an algebraic variety $X$ over $\CC$
is said to be {\em strongly $C^{\infty}$ at $x \in X$} if
there are an open neighborhood $U$ of $x$,
a complex manifold $V$, and a $C^{\infty}$-function $\Phi$ on $V$
such that $U$ is a closed analytic subset of $V$ and
$\phi = \rest{\Phi}{U}$.
If $\phi$ is strongly $C^{\infty}$ at any points of $X$,
then $\phi$ is said to be {\em strongly $C^{\infty}$} on $X$.
We say a continuous hermitian line bundle $\overline{L} = (L, \Vert\cdot\Vert)$
is {\em strongly $C^{\infty}$} if,
for any local basis $l$ of $L$ on any open set, 
$\Vert l \Vert$ is strongly $C^{\infty}$ around there.

\bigskip
First, let us consider the following claim.

\begin{Claim}
There is a commutative diagram:
\[
 \begin{CD}
 X       @<{v}<< X' \\
 @V{f}VV         @VV{f'}V \\
 S       @<<{u}< S'
 \end{CD}
\]
with the following properties.
\begin{enumerate}
\renewcommand{\labelenumi}{(\arabic{enumi})}
\item
$f' : X' \to S'$ is a flat and projective morphism of
algebraic varieties over $\CC$, and $S'$ is non-singular.

\item
$u$ and $v$ are projective birational morphisms.

\item
$v^*(\overline{L}_1), \ldots, v^*(\overline{L}_n)$ and
$v^*(\phi)$ are strongly $C^{\infty}$.
\end{enumerate}
\end{Claim}

Let $\mu : X_1 \to X$ be a resolution of singularities of $X$.
Note that we may take $\mu$ as a projective morphism.
Let $f_1 : X_1 \to Y$ be the composition of morphisms
$X_1 \overset{\mu}{\longrightarrow} X \overset{f}{\longrightarrow} S$.
Then, by applying Raynaud's flattening theorem \cite{Raflat}
to $f_1 : X_1 \to S$, we have the following commutative
diagram:
\[
 \begin{CD}
 X @<{\mu}<< X_1       @<{v_1}<< X' \\
 @V{f}VV @V{f_1}VV         @VV{f'}V \\
 S @= S       @<<{u}< S'
 \end{CD}
\]
where $u$ is a projective birational morphism of algebraic
varieties over $\CC$,
$f' : X' \to S'$ is a flat and projective morphism
of algebraic schemes over $\CC$, and
$X'$ is the main part of $X_1 \times_S S'$.
If $S'$ is not non-singular, taking a desingularization $S'' \to S'$ of $S'$
and replacing $X'$ by $X' \times_{S'} S''$,
we may assume that $S'$ is non-singular.
Since $X$ is integral, so is the generic fiber of $f : X \to S$.
Thus, by virtue of Lemma~\ref{lem:criterion:integral:domain:flat}, 
we can see that $X'$ is an algebraic  variety over $\CC$.
In particular, $v_1$ is birational.
Let $v$ be the composition of morphisms
$X' \overset{v_1}{\longrightarrow} X_1 \overset{\mu}{\longrightarrow} X$.
Then, we have our desired commutative diagram satisfying (1) and (2).
Moreover, since $X_1$ is non-singular, 
$\mu^*(\overline{L}_1), \ldots, \mu^*(\overline{L}_n)$ and
$\mu^*(\phi)$ are strongly $C^{\infty}$.
Thus, so are
$v^*(\overline{L}_1), \ldots, v^*(\overline{L}_n)$ and
$v^*(\phi)$.

\medskip
Let us go back to the proof of our theorem.
The commutative diagram in the above claim
gives rise to the following commutative diagram:
\[
 \begin{CD}
 X       @<{p_1}<< X \times_{S} S' @<{\pi}<< X' \\
 @V{f}VV         @VV{p_2}V @VV{f'}V \\
 S       @<<{u}< S' @= S'
 \end{CD}
\]
where $p_i$ $(i=1, 2)$ is the projection to the $i$-the factor,
and $v = p_1 \cdot \pi$.
Note that by Lemma~\ref{lem:criterion:integral:domain:flat},
$X \times_{S} S'$ is integral.
Since the fiber $(X \times_{S} S')_{s'}$ over $s'$
is canonically isomorphic to 
the fiber $X_{u(s')}$ over $u(s')$ for each $s' \in S'(\CC)$,
we can see
\[
\int_{(X \times_S S')_{s'}} p_1^*(\phi) \bigwedge_{i=1}^n
c_1(p_1^*(\overline{L}_i))  =
 \int_{X_{u(s')}} \phi \bigwedge_{i=1}^n c_1(\overline{L}_i),
\]
which shows us
\addtocounter{Claim}{1}
\begin{equation}
\label{eqn:thm:cont:fiber:integral:1}
\int_{X \times_S S'/S'} p_1^*(\phi) \bigwedge_{i=1}^n
c_1(p_1^*(\overline{L}_i)) =
 u^*\left(
 \int_{X/S} \phi \bigwedge_{i=1}^n c_1(\overline{L}_i)
 \right).
\end{equation}
On the other hand, we would like to see
\addtocounter{Claim}{1}
\begin{equation}
\label{eqn:thm:cont:fiber:integral:2}
 \int_{X'_{s'}}v^*(\phi) \bigwedge_{i=1}^n
c_1(v^*(\overline{L}_i)) =
\int_{(X \times_S S')_{s'}} p_1^*(\phi) \bigwedge_{i=1}^n
c_1(p_1^*(\overline{L}_i))
\end{equation}
for all $s' \in S'(\CC)$.
For this purpose, 
considering a general curve passing through $s'$, we may assume that
$\dim S = 1$.
Then, $X'_{s'}$ and $(X \times_S S')_{s'}$ are Cartier divisors and
$\pi^*((X \times_S S')_{s'}) = X'_{s'}$.
Thus, $\pi_*(X'_{s'}) = (X \times_S S')_{s'}$ as cycles
by virtue of the projection formula and the fact that $\pi$ is birational.
Here we recall the projection formula of integral version:

\begin{Formula}
\label{formula:projection:integral}
Let $g : V \to T$ be a surjective morphism of complete algebraic
varieties over $\CC$,
$\overline{L}_1, \ldots, \overline{L}_{\dim V}$ $C^{\infty}$-hermitian line
bundles on $T$, and let $\phi$ be a $C^{\infty}$-function on $T$.
Then, we have
\[
 \int_V g^*(\phi) \bigwedge_{i=1}^{\dim V} c_1(g^*(\overline{L}_i)) =
 \begin{cases}
  {\displaystyle [\CC(V) : \CC(T)]
 \int_T \phi \bigwedge_{i=1}^{\dim V} c_1(\overline{L}_i)} 
   & \text{if $\dim V = \dim T$}, \\
 0 & \text{if $\dim V > \dim T$}.
 \end{cases}
\]
\end{Formula}

\noindent
Hence, by using the above projection formula,
\[
\int_{X'_{s'}}v^*(\phi) \bigwedge_{i=1}^n
c_1(v^*(\overline{L}_i))  =
\int_{\pi_*(X'_{s'})} p_1^*(\phi) \bigwedge_{i=1}^n
c_1(p_1^*(\overline{L}_i)) =
\int_{(X \times_S S')_{s'}} p_1^*(\phi) \bigwedge_{i=1}^n
c_1(p_1^*(\overline{L}_i)).
\]
Therefore, gathering \eqref{eqn:thm:cont:fiber:integral:1} and
\eqref{eqn:thm:cont:fiber:integral:2}, we get
\[
 \int_{X'/S'}v^*(\phi) \bigwedge_{i=1}^n
c_1(v^*(\overline{L}_i)) 
=u^*\left(
 \int_{X/S} \phi \bigwedge_{i=1}^n c_1(\overline{L}_i)
 \right).
\]
The continuity of the fiber integral of
$C^{\infty}$-forms in strong sense is well known (cf. \cite[Theorem~3.8]{St} and
\cite[Theorem~3.3.2]{King}).
Thus,
\[
 \int_{X'/S'}v^*(\phi) \bigwedge_{i=1}^n
c_1(v^*(\overline{L}_i)) 
\]
is continuous.
Here $u$ is projective and birational.
In particular, $u$ is continuous, closed and surjective.
Hence, we can see that
\[
 \int_{X/S} \phi \bigwedge_{i=1}^n c_1(\overline{L}_i) 
\]
is continuous.
\QED

\begin{comment}
\begin{Lemma}
\label{lem:projection:integral}
Let $g : V \to T$ be a surjective morphism of complete varieties over $\CC$.
Let $\overline{L}_1, \ldots, \overline{L}_n$ be $C^{\infty}$-hermitian line
bundles on $T$, and let $\phi$ be a $C^{\infty}$-function on $T$,
where $n = \dim V$.
Then, we have
\[
 \int_V g^*(\phi) \bigwedge_{i=1}^n c_1(g^*(\overline{L}_i)) =
 \begin{cases}
  {\displaystyle [\CC(V) : \CC(T)]
 \int_T \phi \bigwedge_{i=1}^n c_1(\overline{L}_i)} 
   & \text{if $\dim V = \dim T$}, \\
 0 & \text{if $\dim V > \dim T$}.
 \end{cases}
\]
\end{Lemma}

\Proof
If $V$ and $T$ are non-singular,
this is a well known fact as the projection formula. Otherwise,
let us consider the following commutative diagram:
\[
 \begin{CD}
 V @<{\pi}<< V' \\
 @V{g}VV @VV{g'}V \\
 T @<<{\mu}< T'
 \end{CD}
\]
where $V'$ and $T'$ are non-singular complete varieties over $\CC$,
and $\pi$ and $\mu$ are birational.
Then, by the definition of the integral of $C^{\infty}$-forms on
singular varieties, we can see
\[
 \int_V
 g^*(\phi) \bigwedge_{i=1}^n c_1(g^*(\overline{L}_i))  =
 \int_{V'} \pi^*(g^*(\phi)) \bigwedge_{i=1}^n c_1(\pi^*(g^*(\overline{L}_i))) 
\]
and
\[
 \int_T \phi \bigwedge_{i=1}^n c_1(\overline{L}_i)  =
 \int_{T'} \mu^*(\phi) \bigwedge_{i=1}^n c_1(\mu^*(\overline{L}_i)).
\]
Thus, we get our lemma because
$\pi^*(g^*(\phi)) = {g'}^*(\mu^*(\phi))$ and
$\pi^*(g^*(\overline{L}_i)) = {g'}^*(\mu^*(\overline{L}_i))$
for $i=1, \ldots, n$.
\QED
\end{comment}

\section{Proof of Theorem~\ref{thm:cont:metric:Deligne:pairing}}
\renewcommand{\theequation}{\arabic{section}.\arabic{Claim}}
In this section, we will give the proof
of Theorem~\ref{thm:cont:metric:Deligne:pairing}.
We divide the proof of Theorem~\ref{thm:cont:metric:Deligne:pairing} into three steps.

\bigskip
{\bf Step 1.}\quad
Let $\overline{L}$ be a $C^{\infty}$-hermitian
line bundle on $X$.
If $L$ is $f$-ample, then the metric of
Deligne's pairing $\langle \overline{L} ^{\cdot (n+1)} \rangle(X/S)$ is
continuous.

\medskip
Since
$\langle {\overline{L}^{\otimes m}}^{\cdot (n+1)} \rangle(X/S)
= \langle \overline{L}^{\cdot (n+1)} \rangle(X/S)^{\otimes m^{n+1}}$,
we may assume that $L$ is $f$-very ample. Then,
we have an embedding $\phi : X \hookrightarrow \PP(f_*(L))$ over $S$
with $\phi^*(\OO_{\PP(f_*(L))}(1)) = L$.
Since our question is a local problem with respect to $S$,
we may assume that $f_*(L) \simeq \OO_S^{\oplus N+1}$ for some
$N$. Thus, the above embedding gives rise to
$\phi : X \hookrightarrow \PP^N \times S$ over $S$.
Here we give $L$ a new metric $\normabb'$ arising from
the standard Fubini-Study metric of $\OO_{\PP^N}(1)$.
We denote this new $C^{\infty}$-hermitian line bundle by $\overline{L}'$.

We will show that the metric of
$\langle \overline{L}^{\cdot a}, {\overline{L}'}^{\cdot n+1 - a} \rangle(X/S)$
is continuous for every $0 \leq a \leq n+1$ by induction on $a$.
If $a = 0$, then this holds by virtue of
\cite[Theorem~1.4, Theorem~1.6 and Theorem~3.6]{ZhHRSV}.
Denoting the metric of $\overline{L}$ by $\Vert \cdot \Vert$,
we set
\[
u = \log \left( \frac{\normabb'}{\normabb} \right).
\]
Then, $u$ is a $C^{\infty}$-function on $X$.
On the other hand,
\begin{multline*}
\langle \overline{L}^{\cdot a}, {\overline{L}'}^{\cdot n+1 -a} \rangle(X/S)
\otimes 
\langle \overline{L}^{\cdot a+1}, {\overline{L}'}^{\cdot n-a} 
\rangle(X/S)^{\otimes -1} \\
= \langle \overline{L}^{\cdot a}, \overline{L}' \otimes
\overline{L}^{\otimes -1}, {\overline{L}'}^{\cdot n-a} \rangle(X/S)
= \langle \overline{L}^{\cdot a}, (\OO_X, \exp(u)\vert \cdot \vert_X),
{\overline{L}'}^{\cdot n-a} \rangle(X/S) \\
= \left(\OO_S, \exp
\left(\int_{X/S} u c_1(\overline{L})^{\wedge a} \wedge 
c_1(\overline{L}')^{\wedge n-a}\right)\vert \cdot \vert_S
\right),
\end{multline*}
where for an algebraic variety $T$ over $\CC$,
we denote the canonical metric of $\OO_T$ by $\vert \cdot \vert_T$. 
It follows from Theorem~\ref{thm:cont:fiber:integral} that
${\displaystyle \int_{X/S} u c_1(\overline{L})^{\wedge a} \wedge 
c_1(\overline{L}')^{\wedge n-a}}$
is continuous.
Therefore, by hypothesis of induction,
the metric of 
$\langle \overline{L}^{\cdot a+1}, {\overline{L}'}^{\cdot n-a} \rangle(X/S)$
is continuous.

\bigskip
{\bf Step 2.}\quad
Let $\overline{L}_0, \ldots, \overline{L}_{n}$ be $C^{\infty}$-hermitian
line bundles on $X$. If $L_0, \ldots, L_n$ are $f$-ample, then
the metric of Deligne's pairing 
$\langle \overline{L}_0, \ldots, \overline{L}_{n} \rangle(X/S)$ is
continuous.

\medskip
It is not difficult to see that, for variables $X_0, \ldots, X_n$,
\addtocounter{Claim}{1}
\begin{equation}
\label{eqn:step2:proof:main:theorem}
(n+1)! X_0 \cdots X_n = \sum_{I \subset \{0, 1, \ldots, n\}}
(-1)^{n+1-\#(I)}\left( \sum_{i \in I} X_i \right)^{n+1}.
\end{equation}
This formula shows us that
\[
\langle \overline{L}_0, \ldots, \overline{L}_{n} \rangle(X/S)^{\otimes (n+1)!} =
\bigotimes_{I \subset \{0, 1, \ldots, n\}}
\langle \overline{L}_{I}^{\cdot (n+1)}
\rangle(X/S)^{\otimes (-1)^{n+1-\#(I)}},
\]
where $\overline{L}_{I} = \bigotimes_{i \in I}\overline{L}_{i}$.
Therefore, by Step~1, 
we can see that the metric
of $\langle \overline{L}_0, \ldots, \overline{L}_{n} \rangle(X/S)$
is continuous.

\bigskip
{\bf Step 3.}\quad General case.

\medskip
Let $\overline{L}$ be a $C^{\infty}$-hermitian line bundle on $X$
such that $L$ and $L \otimes L_i^{\otimes -1}$'s  are $f$-ample.
We set $\overline{M}_i = \overline{L} \otimes \overline{L}_i^{\otimes -1}$.
Then, 
\[
\langle \overline{L}_0, \ldots, \overline{L}_{n} \rangle(X/S) =
\langle \overline{L} \otimes \overline{M}_0^{\otimes -1}, \ldots,
\overline{L} \otimes \overline{M}_{n}^{\otimes -1}\rangle(X/S).
\]
Thus, using the linearity of Deligne's pairing and Step~2,
we can conclude the proof of our theorem.
\QED
\renewcommand{\theequation}{\arabic{section}.\arabic{Theorem}.\arabic{Claim}}

\section{Proof of Theorem~\ref{thm:hodge:index}}
\label{sec:application}

Before starting the proof of Theorem~\ref{thm:hodge:index},
let us begin with two lemmas.

\begin{Lemma}
\label{lem:Deligne:pairing}
Let $X$ be a smooth projective curve of genus
$g$ over a field $K$. 
Let us consider Deligne's pairing 
\[
\langle \ , \ \rangle : \Pic(X \times J_{X}) \times \Pic(X \times J_{X})
\to \Pic(J_{X})
\]
with respect to the projection $X \times J_{X} \to J_{X}$,
where $J_{X} = \Pic^0(X)$.
Then, we have the following.
\begin{enumerate}
\renewcommand{\labelenumi}{(\arabic{enumi})}
\item
A line bundle $\langle Q, Q \rangle$ does not depend on the choice
of a universal line bundle $Q$ on $X \times J_{X}$.
In other words, if $Q$ and $Q'$ are universal line bundles on
$X \times J_{X}$, then
$\langle Q, Q \rangle \simeq \langle Q', Q' \rangle$.

\item
We assume that there is $c_0 \in \Pic^1(X)$ with 
$(2g-2)c_0 \sim c_1(\omega_X)$.
\rom{(}This holds if we enlarge the base field $K$.\rom{)}
Let us consider an embedding $\phi : X \to J_{X}$
given by $\phi(x) = x - c_0$.
Let $\Theta_X$ be the theta divisor on $J_{X}$ in terms of $\phi$,
namely, $\Theta_X$ is an ample and symmetric divisor
on $J_{X}$ defined by
\[
 \Theta_X = \{ a \in J_{X} \mid 
\text{$a = \phi(x_1) + \cdots + \phi(x_{g-1})$ 
for some $x_1, \ldots, x_{g-1} \in X$} \}.
\]
If $Q$ is a universal line bundle, then
$\langle Q, Q \rangle \simeq \OO_{J_{X}}(-2 \Theta_X)$.
\end{enumerate}
\end{Lemma}

\Proof
(1) Let $p_1 : X \times J_{X} \to X$ (resp. $p_2 : X \times J_{X} \to J_{X}$)
be the projection onto the first (resp. second) factor.
Then, there is a line bundle $M$ on $J_{X}$ with $Q' = Q \otimes p_2^*(M)$.
Then,
\[
\langle Q', Q' \rangle = \langle Q, Q \rangle \otimes \langle Q, p_2^*(M) \rangle^{\otimes 2}
\otimes \langle p_2^*(M), p_2^*(M) \rangle.
\]
Here $\langle Q, p_2^*(M) \rangle  = \langle p_2^*(M), p_2^*(M) \rangle = \OO_{J_{X}}$
because $Q$ has the degree $0$ along the fibers of $p_2$.
Thus, we get (1).

\medskip
(2) We set
\[
Q_X = \OO_{X \times J_{X}}\left(p_1^*\phi^*(\Theta_X) + p_2^*(\Theta_X) - 
s^*(\Theta_X)\right),  
\] 
where $s : X \times J_{X} \to J_{X}$ is a morphism given by
$s(x, a) = \phi(x) + a$. Then, $Q_X$ is a universal line bundle on
$X \times J_{X}$.

First, we claim that
$\langle Q_X, p_1^*(\omega_X) \rangle = \OO_{J_{X}}$.
It is easy to see that
\[
\langle p_1^*\phi^*(\OO_{J_{X}}(\Theta_X)), p_1^*(\omega_X) \rangle =
\OO_{J_{X}}
\quad\text{and}\quad
\langle p_2^*(\OO_{J_{X}}(\Theta_X)), p_1^*(\omega_X) \rangle =
\OO_{J_{X}}((2g-2)\Theta_X)
\]
On the other hand,
if we set $c_1(\omega_X) \sim \sum_{i=1}^{2g-2} x_i$, then
\[
 \langle s^*(\OO_{J_{X}}(\Theta_X)), p_1^*(\omega_X) \rangle =
\OO_{J_{X}}\left(\sum_{i=1}^{2g-2} T_{\phi(x_i)}^*(\Theta_X)\right),
\]
where, for $a \in J_{X}$, $T_a : J_{X} \to J_{X}$ is a morphism
given by $T_a(z) = z + a$.
Thus, since $\phi(x_1) + \cdots + \phi(x_g) = 0$,
\begin{align*}
\langle Q_X, p_1^*(\omega_X) \rangle & =
\OO_{J_{X}}\left( \sum_{i=1}^{2g-2}
(\Theta_X - T_{\phi(x_i)}^*(\Theta_X))\right) \\
& = \OO_{J_{X}}\left(
\Theta_X - T_{\phi(x_1) + \cdots + \phi(x_{2g-2})}^*(\Theta_X)\right) \\
& = \OO_{J_{X}}.
\end{align*}
Therefore, we get our claim.

Let us start the proof of (2).
First of all, it is well known that
$\det R(p_2)_*(Q_X) = \OO_{J_{X}}(-\Theta_X)$.
On the other hand, by Riemann-Roch theorem,
\[
\det R(p_2)_*(Q_X)^{\otimes 2} = \langle Q_X, Q_X \rangle \otimes
\langle Q_X, p_2^*(\omega_X) \rangle^{\otimes -1} \otimes 
\det R (p_2)_*(\OO_{X \times J_{X}})^{\otimes 2}.  
\]
Here $\det R (p_2)_*(\OO_{X \times J_{X}}) = 
\langle Q_X, p_2^*(\omega_X) \rangle = \OO_{J_{X}}$.
Thus, 
$\langle Q_X, Q_X \rangle = \OO_{J_{X}}(-2 \Theta_X)$.
Therefore, by (1), we have (2).
\QED

\begin{Lemma}
\label{lem:criterion:integral:domain:flat}
Let $A \subseteq R$ be commutative rings such that
$R$ is flat over $A$.
If $A$ and the localization $R_S$ with respect to
a multiplicative subset
$S$ of $A \setminus \{ 0 \}$ are integral domains, then
$R$ is also an integral domain.
\end{Lemma}

\Proof
It is sufficient to show
that the natural homomorphism $R \to R_S$ is injective.
Pick up an element $a \in R$ with $a = 0$ in $R_S$.
Then, there is $s \in S$ with $s a = 0$ in $R$.
Since $s \not= 0$ and $A$ is an integral domain,
$A \overset{\times s}{\longrightarrow} A$ is injective.
Thus, so is $R \overset{\times s}{\longrightarrow} R$
because $R$ is flat over $A$.
Therefore, $a = 0$ in $R$.
\QED

\bigskip
Let us start the proof of Theorem~\ref{thm:hodge:index}.
Fix a universal line bundle $Q$ on $X \times J_X$.
Then, we can find a morphism $q: \mathcal{Y} \to \mathcal{J}$ of
projective arithmetic varieties over $B$ and a $C^{\infty}$-hermitian line
bundle $\overline{\mathcal{Q}}$ on $\mathcal{Y}$
such that $q : \mathcal{Y} \to \mathcal{J}$ coincides with
the projection $X \times J_X \to J_X$ over $K$, and that
$\mathcal{Q}$ is equal to $Q$ over $K$.
If $q$ is not flat, taking Raynaud's flattening theorem \cite{Raflat},
we may assume that $q$ is flat.
Note that after flattening, $\mathcal{Y}$ is still integral
by Lemma~\ref{lem:criterion:integral:domain:flat}.
Let $B_m$ be the closure of the point $[L^{\otimes m}] \in J_X$ in $\mathcal{J}$.
We set
\[
\mathcal{X}_m = q^{-1}(B_m)\quad\text{and}\quad
\overline{\mathcal{L}}_m = \left(\rest{\overline{\mathcal{Q}}}{\mathcal{X}_m}\right)^{\otimes 1/m}
\]
as an element of $\aPic(\mathcal{X}_m) \otimes \QQ$.
Then, since
$\mathcal{X}_m$ is integral by Lemma~\ref{lem:criterion:integral:domain:flat},
we can see that $(\mathcal{X}_m, \overline{\mathcal{L}}_m)$ is a $C^{\infty}$-model of $(X, L)$.
Moreover,
$\mathcal{X}_m \to B_m$ is flat, and $\pi_m : B_m \to B$ is birational.
Here, since Deligne's pairing is compatible with base changes,
\[
\langle \overline{\mathcal{L}}_m^{\otimes m}, 
\overline{\mathcal{L}}_m^{\otimes m} \rangle(\mathcal{X}_m/B_m)
= 
\rest{\langle \overline{\mathcal{Q}}, \overline{\mathcal{Q}} \rangle(\mathcal{Y}/\mathcal{J})}{B_m}.
\]
Thus,
\[
\adeg \left(
\acherncl_1 \left(\langle  \overline{\mathcal{L}}_m^{\otimes m}, 
\overline{\mathcal{L}}_m^{\otimes m}
\rangle(\mathcal{X}_m/B_m)\right)
\cdot \acherncl_1(\pi_m^*(\overline{H}))^d \right) = 
h^{\overline{B}}_{\langle \overline{\mathcal{Q}}, \overline{\mathcal{Q}}
\rangle(\mathcal{Y}/\mathcal{J})}([L^{\otimes m}]).
\]
Clearly,
\[
\adeg \left(
\acherncl_1 \left(\langle  \overline{\mathcal{L}}_m^{\otimes m}, 
\overline{\mathcal{L}}_m^{\otimes m} 
\rangle(\mathcal{X}_m/B_m)\right)
\cdot \acherncl_1(\pi_m^*(\overline{H}))^d \right) = m^2
\adeg \left(
\acherncl_1 (\overline{\mathcal{L}}_m) \cdot \acherncl_1(\overline{\mathcal{L}}_m)
\cdot \acherncl_1(f_m^*(\overline{H}))^d \right),
\]
where $f_m$ is the composition of morphisms $\mathcal{X}_m \to B_m \to B$.
Hence, we get
\[
\adeg \left(
\acherncl_1 (\overline{\mathcal{L}}_m) \cdot \acherncl_1(\overline{\mathcal{L}}_m)
\cdot \acherncl_1(f_m^*(\overline{H}))^d \right) =
\frac{1}{m^2}
h^{\overline{B}}_{\langle \overline{\mathcal{Q}}, \overline{\mathcal{Q}}
\rangle(\mathcal{Y}/\mathcal{J})}([L^{\otimes m}]).
\]
On the other hand, 
the metric of $\langle \overline{\mathcal{Q}}, \overline{\mathcal{Q}}
\rangle(\mathcal{Y}/\mathcal{J})$ is continuous by Theorem~\ref{thm:cont:metric:Deligne:pairing}.
Thus, by \cite[Corollary~3.3.5]{MoArht} and Lemma~\ref{lem:Deligne:pairing},
we can see that
\[
\lim_{n\to\infty}
\frac{1}{m^2} h^{\overline{B}}_{\langle \overline{\mathcal{Q}}, \overline{\mathcal{Q}}
\rangle(\mathcal{Y}/\mathcal{J})}([L^{\otimes m}]) = \hat{h}^{\overline{B}}_{-2\Theta_X}([L]) =
-2 \hat{h}^{\overline{B}}_{\Theta_X}([L]).
\]
Therefore, we get our theorem.

\end{document}